\documentclass[reqno]{lms}

\usepackage{amsmath}
\usepackage{amsfonts}
\usepackage{amssymb,enumerate}
\usepackage[all]{xy}
\usepackage{rotating}
\usepackage{hyperref, color}
\newtheorem{lem}{Lemma}[section]

\newtheorem{prop}[lem]{Proposition}
\newtheorem{thm}[lem]{Theorem}

\newtheorem{intthm}{Theorem}

\newnumbered{defn}[lem]{Definition}
\newnumbered{ex}[lem]{Example}
\newnumbered{question}[lem]{Question}
\newnumbered{problem}[lem]{Problem}
\newnumbered{disc}[lem]{Remark}
\newnumbered{construction}[lem]{Construction}
\newnumbered{notn}[lem]{Notation}
\newnumbered{fact}[lem]{Fact}
\newnumbered{para}[lem]{}
\newnumbered{exer}[lem]{Exercise}
\newnumbered{remarkdefinition}[lem]{Remark/Definition}
\newnumbered{notation}[lem]{Notation}
\newnumbered{step}{Step}
\newnumbered{claim}{Claim}
\newnumbered{convention}[lem]{Convention}
\newnumbered{outline}[lem]{Outline}

\newcommand{\catd}{\mathcal{D}}

\newcommand{\id}{\operatorname{id}}

\newcommand{\depth}{\operatorname{depth}}	
\newcommand{\rank}{\operatorname{rank}}	
\newcommand{\amp}{\operatorname{amp}}

\newcommand{\mspec}{\operatorname{m-Spec}}

\newcommand{\HH}{\operatorname{H}}

\newcommand{\spec}{\operatorname{Spec}}
\newcommand{\s}{\mathfrak{S}}

\newcommand{\im}{\operatorname{Im}}
\newcommand{\shift}{\mathsf{\Sigma}}

\newcommand{\Pic}{\operatorname{Pic}}

\newcommand{\End}{\operatorname{End}}

\newcommand{\ideal}[1]{\mathfrak{#1}}
\newcommand{\m}{\ideal{m}}

\newcommand{\p}{\ideal{p}}

\newcommand{\ol}{\overline}
\newcommand{\wti}{\widetilde}

\newcommand{\supp}{\operatorname{Supp}}

\newcommand{\bbz}{\mathbb{Z}}

\newcommand{\xra}{\xrightarrow}

\newcommand{\into}{\hookrightarrow}

\newcommand{\vf}{\varphi}

\newcommand{\y}{\mathbf{y}}

\newcommand{\x}{\mathbf{x}}

\renewcommand{\geq}{\geqslant}
\renewcommand{\leq}{\leqslant}

\newcommand{\Ext}[4][R]{\operatorname{Ext}_{#1}^{#2}(#3,#4)}	
	
\newcommand{\Rhom}[3][R]{\mathbf{R}\!\operatorname{Hom}_{#1}(#2,#3)}	
\newcommand{\Lotimes}[3][R]{#2\otimes^{\mathbf{L}}_{#1}#3}
\newcommand{\Otimes}[3][R]{#2\otimes_{#1}#3}
\newcommand{\Hom}[3][R]{\operatorname{Hom}_{#1}(#2,#3)}	
\newcommand{\grHom}[3][R]{\operatorname{gr-Hom}_{#1}(#2,#3)}	
\newcommand{\Tor}[4][R]{\operatorname{Tor}^{#1}_{#2}(#3,#4)}

\newcommand{\ba}{\mathbf{a}}

\newcommand{\und}[1]{#1^{\natural}}

\newcommand{\HomA}[2]{\operatorname{Hom}_{A}(#1,#2)}

\newcommand{\RhomA}[2]{\mathbf{R}\!\operatorname{Hom}_{A}(#1,#2)}

\newcommand{\od}{\operatorname{Mod}}
\newcommand{\grod}{\operatorname{Mod}_{\text{gr}}}

\newcommand{\ggl}{\operatorname{GL}}
\newcommand{\ta}{\operatorname{\textsf{T}}}

\newcommand{\en}{\operatorname{End}}

\newcommand{\yext}[4][R]{\operatorname{YExt}_{#1}^{#2}(#3,#4)}	
	
\newcommand{\yexta}[4][R]{\operatorname{Ext}_{#1}^{#2}(#3,#4)}	
\newcommand{\gryexta}[4][R]{\operatorname{gr-Ext}_{#1}^{#2}(#3,#4)}	
\newcommand{\stab}[1][M]{\operatorname{Stab}(#1)}

\newcommand{\Tamod}{\ta^{\od^U(W)}}
\newcommand{\Tamoda}{\ta^{\od^A(W)}}
\newcommand{\Tagrmod}{\ta^{\grod^V(W)}}
\newcommand{\Tagrmoda}{\ta^{\grod^A(W)}}
\newcommand{\Taglm}{\ta^{\ggl(W)_0\cdot M}}
\newcommand{\Tagglm}{\ta^{\ggl(W)\cdot M}}
\newcommand{\Taglum}{\ta^{\ggl(W)_0\cdot \su M}}
\newcommand{\Taglvm}{\ta^{\ggl(W)_0\cdot \sv M}}
\newcommand{\Tagl}{\ta^{\ggl(W)_0}}
\newcommand{\Taggl}{\ta^{\ggl(W)}}
\newcommand{\Fe}[1]{#1[\epsilon]}

\newcommand{\dpic}{\operatorname{DPic}}
\newcommand{\orbit}[1]{\dpic(R)\cdot #1}
\newcommand{\nGor}{\operatorname{nGor}}

\newcommand{\sv}[1]{{}_V\! #1}
\newcommand{\su}[1]{{}_U\! #1}
\newcommand{\gm}{G_\text{m}}

\numberwithin{equation}{lem}

\begin{document}

\bibliographystyle{amsplain}

\title[Geometric representation theory for {DG} algebras]
{Geometric aspects of representation theory for {DG} algebras: answering a question of Vasconcelos}

\author{Saeed Nasseh and Sean Sather-Wagstaff}

\thanks{Sather-Wagstaff  was supported in part by North Dakota EPSCoR, 
National Science Foundation Grant EPS-0814442,
and NSA grant H98230-13-1-0215.}

\date{\today}

\dedication{To Wolmer V. Vasconcelos}

\classno{
13D02, 
13D09, 
13E10, 
14L30, 
16G30} 

\begin{abstract}
We apply geometric techniques from representation theory to the study of homologically finite differential graded (DG) modules $M$ over a finite dimensional, positively graded, commutative DG algebra $U$. In particular, in this setting we prove a version of a theorem of Voigt by exhibiting an isomorphism between the Yoneda Ext group $\operatorname{YExt}^1_U(M,M)$ and a quotient of tangent spaces coming from an algebraic group action on an algebraic variety. As an application, we answer a question of Vasconcelos from 1974 by showing that a local ring has only finitely many semidualizing complexes up to shift-isomorphism in the derived category $\mathcal{D}(R)$.
\end{abstract}

\maketitle

\tableofcontents

\section{Introduction} \label{sec0}

\begin{convention*}\label{conv110211a}
In this paper, $R$ is a commutative  noetherian
ring with identity.
\end{convention*}

The main result of this paper is concerned with
a question of Vasconcelos about \emph{semidualizing $R$-modules},
i.e.,  finitely generated $R$-modules $C$ such that
$\Hom CC\cong R$ and $\Ext iCC=0$ for $i\geq 1$.
These modules were introduced
by Foxby~\cite{foxby:gmarm}
as 
``PG modules of rank 1'', providing
a common generalization of Grothendieck's
canonical modules over Cohen-Macaulay rings~\cite{hartshorne:lc} and the 
projective modules of rank 1. (The P in PG is for ``projective'' and the G is for ``Gorenstein'', referring to Sharp's
Gorenstein modules~\cite{sharp:gmccmlr}. However, the definition makes no assumption about the projective
or injective dimension of $C$.)
Foxby's objective was to provide a single framework that explained useful functorial properties exhibited by
these two seemingly  different classes of modules.

Vasconcelos~\cite{vasconcelos:dtmc} discovered
these objects independently 
(calling them ``spherical modules'') in his investigation of divisors.
Later Golod~\cite{golod:gdagpi} singled them out (also independently, but calling them ``suitable modules''\footnote{
According to private communication with Foxby, the Russian term originally used by Golod has several translations,
including ``comfortable modules''.}) to give a general duality for 
Auslander and Bridger's G-dimension~\cite{auslander:smt}. 
Wakamatsu~\cite{wakamatsu:mtse} similarly identified them in as a useful generalization of tilting modules in representation theory,
calling them, appropriately, 
``generalized tilting modules''. The point in much of this work is that these generalizations give enough additional flexibility
to allow one to verify properties for them that specialize to interesting results about
their progenitors.
For instance, Sather-Wagstaff~\cite{sather:bnsc} was able to use them to get some leverage
on a question of Huneke about Betti numbers of canonical modules, i.e., Bass numbers of  rings.

Avramov and Foxby~\cite{avramov:rhafgd} used  a generalized form of Golod's construction
in their study of G-dimensions of local ring homomorphisms $R\to S$. In particular, they
establish a sweeping extension
of the classical result characterizing the Gorenstein property of $S$ in terms of the Gorensteinness of $R$ and of the closed fibre
of the homomorphism when $S$ is flat over $R$.
In this work, they posed a question that has dogged many of us in the area: 
if two local ring homomorphisms $R\to S\to T$ have finite G-dimension, must the same be true of the
composite map $R\to T$?
This question, with its relation to semidualizing modules (and more generally, semidualizing complexes;
see Appendix~\ref{sec110211a} for background information on these) 
is one of our main reasons for investigating these objects, in addition to the connection with Bass numbers.

In spite of the utility of the semidualizing modules, we know surprisingly little about them.
For instance, we are only now in a position to answer the  following 
aforementioned question of Vasconcelos.

\begin{question}[\protect{\cite[p.~97]{vasconcelos:dtmc}}]\label{q111227a}
If $R$ is local and Cohen-Macaulay, must the set of isomorphism classes
of semidualizing $R$-modules  be finite?
\end{question}

The main result of this paper, stated next, provides a complete answer to Vasconcelos' question.
Note that it does not assume that $R$ is Cohen-Macaulay, and deals with the more general semidualizing
complexes.

\begin{intthm}\label{cor111115a}
If $R$ is a local ring,
then the set  of shift-isomorphism classes of semidualizing
$R$-complexes in 
the derived category
$\catd(R)$  is finite.
\end{intthm}

Section~\ref{sec111227b} is devoted to the proof of this result (contained in~\ref{proof111227c})
and some consequences, including a version of this
result for semi-local rings
(in Theorem~\ref{thm111212a}).

Our proof is inspired by
Christensen and Sather-Wagstaff's~\cite{christensen:cmafsdm} treatment of the special case where $R$ is Cohen-Macaulay and contains a field, so we summarize their proof here.
First,
they reduce to the case where $R$ is complete with
algebraically closed residue field $F$.
They replace $R$ with the finite dimensional $F$-algebra $R/(\mathbf x)$, where $\mathbf x$ is a maximal $R$-sequence.
The desired result then follows from a 
theorem of Happel~\cite{happel:sm}: in this context, there are only finitely many
isomorphism classes of
$R$-modules $C$ of a given length $r$ such that $\Ext 1CC=0$.

Happel's result is proved using geometric techniques.
One parametrizes all  $R$-modules of length $r$  
by an algebraic variety $\od_r^R$
that is acted on by the general linear group
$\ggl_r^F$. The isomorphism class of a module $M$
is  the orbit $\ggl_r^F\cdot M$,
and the tangent space $\ta^{\ggl_r^F \cdot M}_M$
to the orbit $\ggl_r^F \cdot M$ at $M$
is naturally identified with
a subspace of the tangent space $\ta^{\od_r^R}_M$. 
A~theorem of Voigt~\cite{voigt:idteag}
(see also Gabriel~\cite{gabriel:frto} 
and Theorem~\ref{Voigt} below)
provides an isomorphism between $\Ext 1MM$ and the quotient $\ta^{\od_r^R}_M/\ta^{\ggl_r^F \cdot M}_M$.
If $\Ext 1MM=0$, then the orbit $\ggl_r^F \cdot M$ is open in $\od_r^R$.
As $\od_r^R$ is quasi-compact, it can only have finitely many open orbits,
and Happel's result follows.

The idea behind our proof is the same as that of~\cite{christensen:cmafsdm},
with  important differences. First, without the Cohen-Macaulay assumption, given a maximal
$R$-regular sequence $\x$, the quotient $R/(\x)$ is not artinian;
and if $R$ does not contain a field, then there is no reason to expect that the quotient will either.
One cannot remedy this problem by replacing the regular sequence with a system of parameters 
$\y$ (or in fact any sequence generating an ideal of finite colength)
because 
in general
one loses too much homological information in the passage from $R$ to $R/(\y)$.

To deal with this, we employ a technique pioneered by Avramov:
instead of replacing $R$ with $R/(\mathbf{x})$, we use the Koszul complex
$K$ on a minimal generating sequence for the maximal ideal of $R$.
More specifically, we replace $R$ with a related finite dimensional \emph{DG algebra}
over an algebraically closed field $F$.
(See Appendix~\ref{sec110211a} for background information on DG algebras and
DG modules.)

To prove  versions of the results of Happel and Voigt, in Section~\ref{sec161119b}
we develop certain geometric aspects of representation theory for DG 
algebras.
In short, we parametrize all finite dimensional DG $U$-modules $M$ with fixed
underlying graded $F$-vector space $W$
by an algebraic variety $\od^U(W)$.
A product  $\ggl(W)_0$ of
general linear groups acts on this variety 
so that the isomorphism class of $M$
is precisely the orbit $\ggl(W)_0\cdot M$. 
Our version  of Voigt's result for this context is the following;
its proof is 
in~\ref{proof161213a}.

\begin{intthm}\label{tangent space and Ext}
Let $U$ be a finite dimensional DG algebra over an algebraically closed field $F$,
and let $W$ be a finite dimensional graded $F$-vector space.
Given an element $M\in\od^U(W)$, the
quotient of tangent spaces $\Tamod_M/\Taglm_M$
is isomorphic
to the Yoneda Ext group $\yext[U]{1}MM$ 
from Definition~\ref{disc110616a}.\footnote{
See Definition~\ref{DGK2'} for the different but equally important $\Ext[U]1MM$.}
\end{intthm}

As a consequence, if $\yext[U] 1MM=0$, then the orbit $\ggl(W)_0\cdot M$ is open in $\od^U(W)$,
and our version of Happel's result follows;
see 
Proposition~\ref{prop111115a}.

Our first proof of Theorem~\ref{tangent space and Ext} mimicked the proof of Voigt's
original result. In contrast, the proof currently given obtains our result as a corollary of 
Theorem~\ref{grVoigt}, which is itself a corollary of Voigt's original result
and is the subject of Section~\ref{sec161119c}.

\subsection*{Acknowledgments}
We are grateful to \texttt{a-fortiori@mathoverflow.net}, L. L. Avramov, A. Bertram, L.~W.~Christensen,
D. Fulghesu, D. Happel, S. Iyengar, P. J\o rgensen, B. Keller,  P. McNamara, and U. Nagel   for helpful discussions about this work.
We are also grateful to the anonymous referees for many thoughtful suggestions,
including the present proof of Theorem~\ref{tangent space and Ext}.

\section{Graded Version of Voigt's Theorem} \label{sec161119c}

In this section we derive a graded version of Voigt's theorem for use in our proof of Theorem~\ref{tangent space and Ext}
in the next section. We state Voigt's original result in Theorem~\ref{Voigt} below, after developing the necessary notation.
Our graded version is a corollary of this; see Theorem~\ref{grVoigt}.
The ideas here are from~\cite{SGA3,brion:rq,gabriel:frto,happel:sm}.

\begin{notation}\label{notn161119c}
Let $F$ be an algebraically closed field, and set $\gm:=F^\times$,
the multiplicative group over $F$.

Let $A$ be
a finite dimensional associative graded $F$-algebra such that $A_i=0$ for $|i|>q$.
(Note that we do not assume 
multiplicative
commutativity for $A$.)
Let $\dim_F(A_i)=c_i$ for  $i=-q,\ldots, q$. Let
$W:=\bigoplus_{i=0}^s W_i$
be a finite dimensional graded $F$-vector space with $r_i:=\dim_F(W_i)$ for  $i=0,\ldots, s$.
We use the convention $W_j=0$ for all $j\notin\{0,\ldots,s\}$.

Let $\en_F(W)$ denote the set of  $F$-linear endomorphisms of $W$:
\begin{equation}\label{eq161120b}
\en_F(W)=\prod_{i=0}^s\prod_{j=-s}^s\Hom[F]{W_i}{W_{i+j}}
\end{equation}
For each $j\in\bbz$, let $\en_F(W)_j$ denote the set of  $F$-linear endomorphisms of $W$ that are homogeneous of degree $j$:
$$\End_F(W)_{j}=\Hom[F]WW_{j}=\prod_{i=0}^s \Hom[F]{W_i}{W_{i+j}}.$$
\end{notation}

\subsection*{Representations of $A$-Module Structures}

A (left) $A$-module structure on $W$ consists of 
a scalar multiplication 
\begin{equation}\label{eq161120a}
\mu\in\Hom[F]{A\otimes_F W}{W}\cong\prod_{i=0}^q\prod_{j=0}^s\prod_{k=-q-s}^{q+s}\Hom[F]{A_i\otimes_FW_j}{W_{i+j+k}}
\end{equation}
that is associative and unital.
A graded $A$-module structure on $W$ consists of 
a scalar multiplication 
$$\mu\in\Hom[F]{A\otimes_F W}{W}_0\cong \prod_{i=0}^q\prod_{j=0}^s \Hom[F]{A_i\otimes_F W_j}{W_{i+j}}$$
that is associative and unital.

\begin{notation}\label{notn161119d}
Let $\od^A(W)$ denote the set of all maps $\mu$
making $W$ into an $A$-module. 
Sometimes we identify a map $\mu\in\od^A(W)$ with the corresponding  $A$-module
$M=(W,\mu)$.
Let $\grod^A(W)\subseteq\od^A(W)$ denote the set of all maps $\mu$
making $W$ into a graded $A$-module. 
\end{notation}

\begin{disc}\label{disc161119b}
The multiplicative group $\gm$ acts on $W$; on homogeneous elements $w$, this action has the form $\alpha\cdot w=\alpha^{|w|}w$
where the right-hand expression uses the scalar multiplication of $F$ on $W$.
One defines an action of $\gm$ on $A$ using the same formula.
The actions of $\gm$ on $W$ and $A$ allow
$\gm$ to act on
$\od^A(W)$ by the formula
$\alpha\mu:=\alpha\circ \mu\circ (\alpha^{-1}\otimes_F \alpha^{-1})$.
Because $F$ is infinite, one checks that $\mu$ respects the grading on $W$ if and only if $\mu$ is fixed by this action,
that is, $\grod^A(W)$ is the set of fixed points for this action:
$\grod^A(W)=\od^A(W)^{\gm}$.
\end{disc}

We next describe geometric structures on the sets $\od^A(W)$ and $\grod^A(W)$.

\begin{disc}\label{disc161119a}
The $F$-vector space 
$\Hom[F]{A\otimes_F W}{W}$
from~\eqref{eq161120a}
has dimension $b:=(\sum_ic_i)(\sum_jr_j)^2$, so a map $\mu$
corresponds to an element of the affine space $\mathbb{A}^{b}_F$. The condition that $\mu$ be  associative and unital
is equivalent to the entries of the matrices
representing $\mu$ satisfying certain fixed polynomials over $F$. Thus, the
set $\od^A(W)$  is a Zariski-closed subset of $\mathbb{A}^{b}_F$.

In the decomposition~\eqref{eq161120a},
the maps $\mu\in\Hom[F]{A\otimes_F W}{W}$ that are homogeneous of degree 0
are exactly the ones where the matrices from $\Hom[F]{A_i\otimes_FW_j}{W_k}$ with $k\neq i+j$ are all 0.
In other words, the space $\Hom[F]{A\otimes_F W}{W}_0$ is a linear subspace of $\Hom[F]{A\otimes_F W}{W}$
of dimension $b':=\sum_{t} \sum_{i} c_ir_{t-i}r_t$. 
We identify the affine space $\mathbb{A}^{b'}_F$ with the corresponding linear subvariety of $\mathbb{A}^{b}_F$;
the defining conditions for $\grod^A(W)$ and $\od^A(W)$,
namely the associative and unital conditions, then say 
that
$\grod^A(W)=\mathbb{A}^{b'}_F\cap\od^A(W)$.
\end{disc}

\begin{notation}\label{notn160311b}
Let $\ggl(W)$ denote the set of  $F$-linear automorphisms of $W$, that is, the invertible elements of $\en_F(W)$.
Let $\ggl(W)_0$ denote the set of  $F$-linear automorphisms of $W$ that are homogeneous of degree 0, that is, the invertible elements of $\en_F(W)_0$; so $\ggl(W)_0=\ggl(W)\cap\en_F(W)_0$.
\end{notation}

\begin{disc}\label{disc160311b}
An element $\alpha\in \ggl(W)$  is an element of 
$\en_F(W)$
with a multiplicative inverse. The vector space $\en_F(W)$ has dimension $\sigma:=(\sum_ir_i)^2$,
so the map $\alpha$ corresponds to an element of the affine space $\mathbb{A}^{\sigma}_F$. 
The invertibility of $\alpha$ is  an open condition, given by the non-vanishing of the determinant.
Thus, the group $\ggl(W)$  is a Zariski-open subset of $\mathbb{A}^{\sigma}_F$;
in particular, $\ggl(W)$ is non-singular  of dimension $\sigma$.

Note that the maps defining the operations in $\ggl(W)$ are regular, that is, 
they are determined by rational functions on the entries of the matrix with non-vanishing denominators on $\ggl(W)$. 
For multiplication, this is from the definition of matrix multiplication in terms of quadratic forms on the entries.
For inversion, this follows from the formula $\alpha^{-1}=\det(\alpha)^{-1}\operatorname{adj}(\alpha)$.
Thus,
$\ggl(W)$ is an algebraic group.

In the decomposition~\eqref{eq161120b},
the maps $\alpha\in\en_F(W)$ that are homogeneous of degree 0
are exactly the ones with a block-diagonal form, that is,
where the matrices from $\Hom[F]{W_i}{W_j}$ with $j\neq i$ are all 0.
In other words, the space $\en_F(W)_0$ is a linear subspace of $\en_F(W)$
of dimension $e:=\sum_i r_i^2$. 
We identify the affine space $\mathbb{A}^{e}_F$ with the corresponding linear subvariety of $\mathbb{A}^{\sigma}_F$.
The defining conditions for $\ggl(W)_0$ and $\ggl(W)$,
namely the invertibility conditions, then say 
that
$\ggl(W)_0=\mathbb{A}^{e}_F\cap\ggl(W)$.
In particular, $\ggl(W)_0$ is non-singular  of dimension $e$.
Moreover, as with $\ggl(W)$, the subset $\ggl(W)_0\subseteq\ggl(W)$ is an algebraic subgroup.
(Alternately, one sees part of this by viewing $\ggl(W)_0$ as the product $\ggl(W_0)\times\cdots\times\ggl(W_s)$,
since each $\ggl(W_i)$ is a non-singular algebraic group.)

Note that the action of $\gm$ on $W$ from Remark~\ref{disc161119b} allows us to identify
$\gm$ naturally with a subgroup of $\ggl(W)_0$, hence of $\ggl(W)$.
Under this identification, 
because $F$ is infinite, one checks that
$\ggl(W)_0$ is the centralizer of $\gm$ in $\ggl(W)$.
\end{disc}

\subsection*{General Linear Group Action on $\od^A(W)$}

Next, we describe an action of the group 
$\ggl(W)$
from Notation~\ref{notn160311b}
on $\od^A(W)$ by conjugation.

\begin{notation}\label{notn161120a}
Let $\alpha\in \ggl(W)$. For every $\mu\in\od^A(W)$, we define 
\begin{equation}\label{eq161120c}
\alpha\cdot\mu:=\alpha\circ \mu\circ (A\otimes_F \alpha^{-1}).
\end{equation}
\end{notation}

Let $\alpha\in \ggl(W)$ and $\mu\in\od^A(W)$.
It is  straightforward to show 
that~\eqref{eq161120c}
describes an $A$-module structure for $W$, so
$\alpha\cdot\mu\in\od^A(W)$.
From the definition of $\alpha\cdot\mu$, it follows readily that this describes a
$\ggl(W)$-action on $\od^A(W)$.
In addition, the maps defining 
this action
are regular, essentially by the second paragraph
of Remark~\ref{disc160311b}.
In the language of~\cite{brion:rq}, this says that we have an algebraic action of $\ggl(W)$ on $\od^A(W)$.
Furthermore, this restricts to an algebraic action of $\ggl(W)_0$ on $\grod^A(W)$:
if $\alpha\in \ggl(W)_0\subseteq\ggl(W)$ and $\mu\in\grod^A(W)\subseteq\od^A(W)$, then 
$\alpha\cdot\mu\in\grod^A(W)$.

We now have the notation needed to state Voigt's original theorem.
Note that we will apply it to the algebra $A$ and the vector space $W$
with the gradings forgotten. 

\begin{thm}[Voigt~\cite{voigt:idteag}]
\label{Voigt}
Let $B$ be a finite dimensional  associative algebra over an algebraically closed field $F$,
and let $X$ be a finite dimensional  $F$-vector space.
Given an element $N\in\od^B(X)$, there is a natural identification of the tangent space $\ta^{\ggl(X)\cdot N}_N$
as a subspace of $\ta^{\od^B(X)}_N$, and under this identification
the
quotient of tangent spaces $\ta^{\od^B(X)}_N/\ta^{\ggl(X)\cdot N}_N$
is isomorphic
to the  Ext group $\yexta[B]{1}NN$.\footnote{In contrast with the DG-setting
of Theorem~\ref{tangent space and Ext}, there is no need for the notation YExt in the non-DG-setting.}
\end{thm}

Next, we describe the orbits of the action of $\ggl(W)$ on $\od^A(W)$.
Let $\alpha\in \ggl(W)$ and $\mu\in\od^A(W)$, and set $\widetilde{\mu}:=\alpha\cdot\mu$.
It is  straightforward to show that  the map $\alpha$ gives an
$A$-module isomorphism $\alpha\colon (W,\mu)\xra\cong(W,\wti\mu)$.
Furthermore, if $\alpha\in \ggl(W)_0$ and $\mu\in\grod^A(W)$, then the isomorphism 
$\alpha\colon (W,\mu)\xra\cong(W,\wti\mu)$ is homogeneous.

Conversely, given another element $\mu'\in\od^A(W)$,
if there is an $A$-module isomorphism $\beta\colon(W,\mu)\xra\cong(W,\mu')$, then
$\beta\in \ggl(W)$ and $\mu'=\beta\cdot\mu$.
Furthermore, if $\mu,\mu'\in\grod^A(W)$ and if $\beta$ is a homogeneous isomorphism,
then
$\beta\in \ggl(W)_0$ and $\mu'=\beta\cdot\mu$.
Thus, we have the following.

\begin{prop}\label{prop161120a}
The orbits in $\od^A(W)$ under the action of $\ggl(W)$
are in bijection with the isomorphism classes
of $A$-module structures on $W$.
The orbits in $\grod^A(W)$ under the action of $\ggl(W)_0$
are in bijection with the isomorphism classes
of graded $A$-module structures on $W$.
\end{prop}

\subsection*{Tangent Spaces}

\begin{notation}\label{notn160312d'}
The Zariski tangent space of $\od^A(W)$ at an element $M=\mu\in\od^A(W)$ is denoted $\Tamoda_M$.
Similarly, for any $\alpha\in\ggl(W)$, we have tangent spaces 
$\Taggl_\alpha$ and $\Tagglm_{\alpha\cdot M}$, and similarly for the homogeneous cases.

Set $F[\epsilon]:=F\epsilon \oplus  F$, where  $\epsilon^2=0$ and $|\epsilon|=0$.
We set 
\begin{gather*}
A[\epsilon]:=\Otimes[F]{F[\epsilon]}{A}\cong A\epsilon \oplus  A\cong A \oplus  A\\
W[\epsilon]:=\Otimes[F]{F[\epsilon]}{W}\cong W\epsilon \oplus W\cong 
W \oplus W.
\end{gather*}
We write elements of $\Fe A$ as
$u\epsilon+v=(\epsilon\otimes u)+(1\otimes v)$, and similarly for $\Fe W$.
\end{notation}

\begin{disc}\label{disc160312a'}
The vector space $A[\epsilon]$ is a graded $F[\epsilon]$-algebra using the  multiplication
$(\zeta\otimes u) (\zeta'\otimes u'):=(\zeta\zeta')\otimes (uu')$.
In other words, since $\epsilon^2=0$,
this works out to be
$$(u\epsilon+v)(u'\epsilon+v')=(uv'+vu')\epsilon+vv'.$$
In particular, $\epsilon u=u\epsilon$ for all $u\in A$.
Similarly, $\Fe W$ is a graded $\Fe F$-module.
\end{disc}

\begin{disc}\label{disc160312b'}
As above,
an  $\Fe A$-module structure on $\Fe W$ consists of a scalar multiplication 
from $\Hom[\Fe F]{\Fe A\otimes_{\Fe F} \Fe W}{\Fe W}$, and a graded module structure
comes from $\Hom[\Fe F]{\Fe A\otimes_{\Fe F} \Fe W}{\Fe W}_0$.

For example, given a (graded) $A$-module structure $M=\mu$ on $W$, 
we have a trivial (graded) $\Fe A$-module structure $\Fe M=\Otimes[F]{\Fe F}{M}$ on $\Fe W$.
In terms of the tensor descriptions of $\Fe A$ and $\Fe M$, this is given as
$(\xi\otimes a)(\zeta\otimes m)=(\xi\zeta)\otimes (am)$.
Writing elements of $\Fe A$ as $u\epsilon+v$ and similarly for $\Fe M$, 
this reads $(u\epsilon+v)(m\epsilon+n)=(un+vm)\epsilon+vn$.
\end{disc}

\begin{notation}\label{notn160312a'}
Let $\od^{\Fe A}(\Fe W)$ denote the set of all maps $\mu$
making $\Fe W$ into an $\Fe A$-module, and
let $\grod^{\Fe A}(\Fe W)\subseteq\od^{\Fe A}(\Fe W)$ denote the subset of all maps $\mu$
making $\Fe W$ into a graded $\Fe A$-module. 
\end{notation}

\begin{disc}\label{disc160312c'}
The natural ring epimorphism $\Fe F\to F$ induces  well-defined maps
$\od^{\Fe A}(\Fe W)\xra\kappa \od^{A}(W)$ 
and 
$\grod^{\Fe A}(\Fe W)\xra{\kappa'} \grod^{A}(W)$
obtained by reducing modulo $\epsilon$.
The point is that 
any unital, associative (graded) scalar multiplication of $\Fe A$ on $\Fe W$ 
reduces modulo $\epsilon$ to 
a unital, associative  (graded)  scalar multiplication of $A$ on $W$.

The tangent space $\Tamoda_M$ of $\od^A(W)$ at an element $M\in\od^A(W)$ is the 
fibre of $\kappa$ over $M$
$$\Tamoda_M=\{N\in\od^{\Fe A}(\Fe W)\mid\kappa(N)=M\}.$$
See, e.g., \cite[VI.1.3]{eisenbud:gs}.
In other words, the tangent space $\Tamoda_M$ is the set of elements of $\od^{\Fe A}(\Fe W)$ that reduce to $M$ modulo $\epsilon$.
A similar conclusion holds for $\Tagrmoda_M$.
\end{disc}

\begin{lem}\label{lem161213a}
Let $M\in\grod^A(W)$.
Then $\Tagrmoda_M$ is naturally a subspace of $\Tamoda_M$.
Moreover, $\gm$ acts on $\Tamoda_M$ by linear maps (that is, by automorphisms) so that
$\Tagrmoda_M=(\Tamoda_M)^{\gm}$.
\end{lem}

\begin{proof}
Remark~\ref{disc161119b} shows how to define an action of $\gm$ on $\od^{\Fe A}(\Fe W)$.
It is straightforward to show that the subspace $\Tamoda_M\subseteq\od^{\Fe A}(\Fe W)$ is invariant under this action. 
One checks readily that $\gm$ acts on $\Tamoda_M$ by linear maps.

The invariance of $\Tamoda_M$ yields the last equality in the next display
\begin{align*}
\Tagrmoda_M
&=\{N\in\grod^{\Fe A}(\Fe W)\mid\kappa'(N)=M\}\\
&=\Tamoda_M\cap\grod^{\Fe A}(\Fe W)\\
&=\Tamoda_M\cap(\od^{\Fe A}(\Fe W))^{\gm}\\
&=(\Tamoda_M)^{\gm}.
\end{align*}
The third equality here is from Remark~\ref{disc161119b}, applied to the algebra $\Fe A$.
The other equalities are from Remark~\ref{disc160312c'}.
\end{proof}

\begin{notation}\label{notn160312c'}
As above,
let $\ggl(\Fe W)$ denote the set of  $\Fe F$-linear automorphisms of $\Fe W$, that is, the 
invertible elements of $\en_{\Fe F}(\Fe W)$,
and let $\ggl(\Fe W)_0$ denote the set of  $\Fe F$-linear automorphisms of $\Fe W$ that are homogeneous of degree 0, that is, the 
invertible elements of $\en_{\Fe F}(\Fe W)_0$.

For each $\alpha\in \ggl(\Fe W)$ and $\mu\in\od^{\Fe A}(\Fe W)$, we define 
$$\alpha\cdot\mu:=\alpha\circ \mu\circ (\Fe A\otimes_{\Fe F} \alpha^{-1}).$$
It follows readily that this 
describes a
$\ggl(\Fe W)$-action on $\od^{\Fe A}(\Fe W)$.
The same formula
describes a
$\ggl(\Fe W)_0$-action on $\grod^{\Fe A}(\Fe W)$.
\end{notation}

The next result is proved like Lemma~\ref{lem161213a}.
The subsequent result, however, requires slightly more effort.

\begin{lem}\label{lem161213b}
The tangent space $\Tagl_{\id_W}$ is naturally a subspace of $\Taggl_{\id_W}$.
Moreover, $\gm$ acts on $\Taggl_{\id_W}$ by linear maps so that
$\Tagl_{\id_W}=(\Taggl_{\id_W})^{\gm}$.
\end{lem}

\begin{lem}\label{lem161213c}
Let $M\in\grod^A(W)$. 
Then the tangent space $\Taglm_M$ is naturally a subspace of $\Tagglm_M$.
Moreover, $\gm$ acts on $\Tagglm_M$ by linear maps so that
$\Taglm_M=(\Tagglm_M)^{\gm}$.
\end{lem}

\begin{proof}
Consider the map $\ggl(W)\xra{\!\cdot M\!}\od^A(W)$
coming from the action of $\ggl(W)$ on $\od^A(W)$, and the induced linear transformation $\Taggl_{\id_W}\xra{D_M}\Tamoda_M$
of tangent spaces.
We describe the  map $D_M$, using the following diagram as a guide.
$$\xymatrix{
\Taggl_{\id_W} \ar[r]^-{D_M}\ar[d]_{\subseteq}
&\Tamoda_M\ar[d]^{\subseteq} \\
\ggl(\Fe W)\ar[r]^-{\cdot\Fe M}&\od^{\Fe A}(\Fe W)
}$$
Identify $\Taggl_{\id_W}$, as above, with the set of elements of $\ggl(\Fe W)$
that reduce to $\id_W$ modulo $\epsilon$. 
The action of $\ggl(\Fe W)$ on $\od^{\Fe A}(\Fe W)$ yields a map
$$\ggl(\Fe W)\xra{\cdot \Fe M}\od^{\Fe A}(\Fe W).$$
The restriction of this map to $\Taggl_{\id_W}$ lands in $\Tamoda_M$, and the restricted map
$\Taggl_{\id_W}\to\Tamoda_M$ is exactly $D_M$. 
In other words, for each $\alpha\in\Taggl_{\id_W}$, we have $D_M(\alpha)=\alpha\cdot\Fe M\in\Tamoda_M$.

Recall that~\cite[Theorem~2.7(ii)]{brion:rq} shows that $\ggl_F(W)\cdot M$ is a locally closed subset of
$\od^A(W)$. Hence, the tangent space $\Tagglm_M$ is naturally identified with
a subspace of $\Tamoda_M$. 
Part of the proof of Voigt's 
Theorem~\ref{Voigt}
states that, under this identification, we have 
\begin{equation}\label{eq161214a}
D_M(\Taggl_{\id_W})=\Tagglm_M.
\end{equation}

From the above description of the map $D_M$, it is straightforward to show that $D_M$ is $\gm$-equivariant.
It follows that the subspace $\Tagglm_M$ is invariant under the action of $\gm$ on $\Tamoda_M$.
Since Lemma~\ref{lem161213b} implies that $\gm$ acts on $\Taggl_{\id_W}$ by linear maps,
we also conclude that $\gm$ acts on the invariant subspace $\Tamoda_M$ by linear maps.

Next, consider
the map $\ggl(W)_0\xra{\!\cdot M\!}\grod^A(W)$
and the induced linear transformation $\Tagl_{\id_W}\xra{D_M}\Tagrmoda_M$.
The above argument translates directly to show
that
$\Taglm_M$ is naturally identified with
a subspace of $\Tagrmoda_M$, and under this identification we have 
\begin{equation}\label{eq161214b}
D_M(\Tagl_{\id_W})=\Taglm_M.
\end{equation}

Recall that the group $\gm$ is linearly reductive. 
It follows that the functor $(-)^{\gm}$ is exact on the category with objects equal to the
finite dimensional $F$-vector spaces with $\gm$-actions by linear maps and with morphisms equal to the equivariant linear maps.
Lemma~\ref{lem161213b} says that $\gm$ acts on $\Taggl_{\id_W}$ by linear maps,
so the exactness of $(-)^{\gm}$ yields
the third equality in the following sequence.
\begin{align*}
\Taglm_M
&=D_M(\Tagl_{\id_W})\\
&=D_M((\Taggl_{\id_W})^{\gm})\\
&=(D_M(\Taggl_{\id_W}))^{\gm}\\
&=(\Tagglm_M)^{\gm}
\end{align*}
The first equality here  
is~\eqref{eq161214b},
the second equality is from 
Lemma~\ref{lem161213b},
and the fourth equality  
is~\eqref{eq161214a}.
\end{proof}

\subsection*{Extensions and our Graded Version of Voigt's Theorem}

\begin{disc}\label{disc161121a}
Let $M=\mu\in\grod^A(W)$. 
The elements of $\yexta[A]1MM$ are equivalence classes of short exact sequences
$$0\to M\to X\to M\to 0.$$
Since the underlying vector space for $M$ is the graded vector space $W$,  such sequences
are split over $F$, 
and
they can be expressed in the form
\begin{equation}
\label{eq161121a}
\xi=\qquad 0\to (W,\mu)\xra i(\Fe W,\nu)\xra p(W,\mu)\to 0
\end{equation}
where $i(w)=w\epsilon$ and $p(w\epsilon+x)=x$.

Within this form, the $F$-vector space structure on $\yexta[A]1MM$ takes the following form.
For an element $\lambda\in F$, we have
\begin{equation*}
\lambda\xi=\qquad 0\to (W,\mu)\xra i(\Fe W,\nu^\lambda)\xra p(W,\mu)\to 0
\end{equation*}
where
$\nu^\lambda(a\otimes(w\epsilon+x)):=\nu(a\otimes(w\epsilon+\lambda x))-\mu(a\otimes(\lambda x))+\mu(a\otimes x)$.
For a second extension
\begin{equation*}
\zeta=\qquad 0\to (W,\mu)\xra i(\Fe W,\tau)\xra p(W,\mu)\to 0
\end{equation*}
the Baer sum in $\yext[A]1MM$ gives
\begin{equation*}
\xi+\zeta=\qquad 0\to (W,\mu)\xra i(\Fe W,\nu\oplus\tau)\xra p(W,\mu)\to 0
\end{equation*}
where
$(\nu\oplus\tau)(a\otimes(w\epsilon+x)):=\nu(a\otimes(w\epsilon+x))+\tau(a\otimes x)-\mu(a\otimes x)$.

Since $M$ is graded, the multiplicative group $\gm$ acts on $\yexta[A]1MM$ as
\begin{equation*}
\alpha \xi=\qquad 0\to (W,\mu)\xra i(\Fe W,\alpha\nu)\xra p(W,\mu)\to 0
\end{equation*}
where the notation $\alpha\nu$ is from Remark~\ref{disc161119b}.
One checks readily that this describes a group action  and furthermore
that $\gm$ acts on $\yexta[A]1MM$ by linear maps.

Moreover, the map $\eta\colon \Tamoda_M\to\yexta[A]1MM$ from Voigt's 
Theorem~\ref{Voigt}
sends an element
$\nu\in\Tamoda_M$ to the extension~\eqref{eq161121a}.
Thus,
it follows  from the definitions that the map $\eta$ is 
$\gm$-equivariant.
\end{disc}

\begin{defn}\label{defn161121a}
Let $M=\mu\in\grod^A(W)$. An extension~\eqref{eq161121a} of $M$ by $M$ is \emph{graded}
when the modules in it are graded as $A$-modules and the maps are homogeneous. 
Let $\gryexta[A]1MM$ denote the set of equivalence classes (under the usual equivalence relation)
of graded extensions of $M$ by $M$.
\end{defn}

\begin{lem}\label{lem161213d}
Let $M=\mu\in\grod^A(W)$. 
Then $\gryexta[A]1MM$ is naturally a subspace of $\yexta[A]1MM$
so that
$\gryexta[A]1MM=\yexta[A]1MM^{\gm}$.
\end{lem}

\begin{proof}
In an extension~\eqref{eq161121a},
the module $M=(W,\mu)$ is  graded by assumption, as are the maps $i$ and $p$.
Thus, such an extension is graded if and only if the module $(\Fe W,\nu)$ is graded.
That is, $\xi\in\gryexta[A]1MM$ if and only if $\nu$ is graded;
by Remark~\ref{disc161119b}, $\nu$ is graded if and only if it is fixed by the action of $\gm$,
and this is so if and only if $\xi$ is fixed by the action of $\gm$.
In summary, this shows that $\gryexta[A]1MM=\yexta[A]1MM^{\gm}$.
\end{proof}

We are finally in a position to prove our graded version of Voigt's Theorem.

\begin{thm}\label{grVoigt}
Let $A$ be a finite dimensional graded associative algebra over an algebraically closed field $F$,
and let $W$ be a finite dimensional graded $F$-vector space.
Given an element $M\in\grod^A(W)$, there is a natural identification of the tangent space $\Taglm_M$
as a subspace of $\Tagrmoda_M$, and under this identification
the
quotient of tangent spaces $\Tagrmoda_M/\Taglm_M$
is isomorphic
to 
$\gryexta[A]{1}MM$.
\end{thm}

\begin{proof}
Voigt's 
Theorem~\ref{Voigt}
provides a short exact sequence
$$0\to \Tagglm_M\xra\subseteq\Tamoda_M\xra\eta\yexta[A]1MM\to 0.$$
The multiplicative group $\gm$ acts on these spaces by linear maps,
and $\eta$ is $\gm$-equivariant;
as we have noted above, the fact that $\gm$ is linearly reductive implies that the induced sequence
$$0\to (\Tagglm_M)^{\gm}\xra\subseteq(\Tamoda_M)^{\gm}\xra{\eta^{\gm}}
\yexta[A]1MM^{\gm}
\to 0$$
is exact;
see Remark~\ref{disc161121a} and
Lemmas~\ref{lem161213a}, \ref{lem161213c}, and~\ref{lem161213d}. These results also
show that
this sequence is
of the form
$$0\to \Taglm_M\xra\subseteq\Tagrmoda_M\to\gryexta[A]1MM\to 0$$
hence the desired conclusions.
\end{proof}

\section{DG Version of Voigt's Theorem} \label{sec161119b}

In this section, we 
prove
Theorem~\ref{tangent space and Ext} from the introduction; 
see~\ref{proof161213a}.
The ideas here 
originate in~\cite{SGA3,brion:rq,gabriel:frto,happel:sm}.
See Appendix~\ref{sec110211a} for background information on DG modules.

\begin{notation}\label{notn160311a}
Let $F$ be an algebraically closed field, and $U$
a finite dimensional  DG $F$-algebra such that $U_i=0$ for $i>q$ and for $i<0$.
Let $\dim_F(U_i)=n_i$ for  $i=0,\ldots, q$. Let
$W:=\bigoplus_{i=0}^s W_i$
be as in Notation~\ref{notn161119c}.
\end{notation}

\subsection*{Representations of DG-Module Structures}

A DG $U$-module structure on $W$ consists of two pieces of data: a differential 
$$\partial\in \End_F(W)_{-1}=\Hom[F]WW_{-1}=\prod_{i=0}^s \Hom[F]{W_i}{W_{i-1}}$$
and a scalar multiplication 
$$\mu\in\Hom[F]{U\otimes_F W}{W}_0\cong \prod_{i=0}^q\prod_{j=0}^s \Hom[F]{U_i\otimes_F W_j}{W_{i+j}}.$$

\begin{notation}\label{notn160311g}
Let $\od^U(W)$ denote the set of all ordered pairs $(\partial,\mu)$
making $W$ into a DG $U$-module. 
Sometimes we identify an ordered pair $(\partial,\mu)\in\od^U(W)$ with the corresponding DG $U$-module
$M=(W,\partial,\mu)$.
\end{notation}

We next describe a geometric structure on the set $\od^U(W)$.

\begin{disc}\label{disc160311a}
As in Section~\ref{sec161119c},
a differential $\partial$ 
on $W$
corresponds to an element of the affine space $\mathbb{A}^{d}_F$ where
$d:=\sum_i r_ir_{i-1}$. The vanishing condition $\partial\partial=0$ is equivalent to the entries of the matrices
representing $\partial$ satisfying certain fixed homogeneous quadratic polynomial equations over $F$. Hence, the
set of all  differentials on $W$ is a Zariski-closed subset of $\mathbb{A}^{d}_F$.
Also,
a map $\mu$
corresponds to an element of the affine space $\mathbb{A}^{d'}_F$ where
$d':=\sum_{c} \sum_{i} n_ir_{c-i}r_c$. The condition that $\mu$ be an associative, unital,  cycle satisfying the Leibniz Rule
is equivalent to the entries of the matrices
representing $\partial$ and $\mu$ satisfying certain fixed polynomials over $F$. Thus, the
set $\od^U(W)$  is a Zariski-closed subset of $\mathbb{A}^{d}_F\times \mathbb{A}^{d'}_F\cong \mathbb{A}^{d+d'}_F$.
\end{disc}

\subsection*{General Linear Group Action on $\od^U(W)$}

Next, we describe an action of the group $\ggl(W)_0$ from Notation~\ref{notn160311b}
on $\od^U(W)$ by conjugation.

\begin{notation}\label{notn160311f}
Let $\alpha\in \ggl(W)_0$. For every $(\partial,\mu)\in\od^U(W)$, we define 
\begin{equation}\label{eq160315a}
\alpha\cdot(\partial,\mu):=(\alpha\circ \partial\circ \alpha^{-1},\alpha\circ \mu\circ (U\otimes_F \alpha^{-1})).
\end{equation}
\end{notation}

Let $\alpha\in \ggl(W)_0$ and $M=(\partial,\mu)\in\od^U(W)$.
As in Section~\ref{sec161119c},
this describes an algebraic action of
$\ggl(W)_0$
on $\od^U(W)$.
Thus,  
most of 
the following result is from~\cite[Proposition~2.7(ii)]{brion:rq}.

\begin{lem}\label{prop160312a}
Let $M\in\od^U(W)$ be given.
\begin{enumerate}[\rm(a)]
\item\label{prop160312a1}
The stabilizer (i.e., the isotropy group)
$$\stab=\{\alpha\in\ggl(W)_0\mid\alpha\cdot M=M\}$$
is closed in $\ggl(W)_0$ and non-singular.
\item\label{prop160312a2}
The orbit 
$$\ggl(W)_0\cdot M=\{\alpha\cdot M\mid\alpha\in\ggl(W)_0\}$$ 
is a locally closed, non-singular subvariety\footnote{In~\cite{brion:rq}, varieties are
not necessarily irreducible.} 
of  $\od^U(W)$.
All connected components of $\ggl(W)_0\cdot M$ have dimension $\dim(\ggl(W)_0)-\dim(\stab)$.
\end{enumerate}
\end{lem}

\begin{proof}
By~\cite[Proposition~2.7]{brion:rq}, we need only show that 
$\stab$ is non-singular.
Since $F$ is algebraically closed, it suffices to show
that $\stab$ is regular. As $\ggl(W)_0$ is regular, to show
that $\stab\subseteq \ggl(W)_0$ is regular it is enough
to
show that $\stab$ is defined by linear equations.
To find these linear equations, note that the stabilizer condition $\alpha\cdot M=M$
is equivalent to the conditions
$\partial=\alpha\circ \partial\circ \alpha^{-1}$ and
$\mu=\alpha\circ \mu\circ (U\otimes_F \alpha^{-1})$,
that is,
$\partial\circ \alpha=\alpha\circ \partial$ and
$\mu\circ (U\otimes_F \alpha)=\alpha\circ \mu$;
since the matrices defining $\partial$ and $\mu$ are fixed,
these equations are described by a system of linear equations
in the variables describing $\alpha$.
Thus, 
$\stab$
is non-singular.
\end{proof}

The next result is verified like Proposition~\ref{prop161120a}.

\begin{prop}\label{prop160311a}
The orbits in $\od^U(W)$ under the action of $\ggl(W)_0$
are in bijection with the isomorphism classes
of DG $U$-module structures on $W$.
\end{prop}

\subsection*{Translating from the DG Setting to the Graded Setting}

The following notation remains in effect for the 
rest of
section.

\begin{notn}\label{notn161119a}
Let $\delta$ be an exterior variable of degree $-1$, so $\delta^2=0$,
and set $V=\und{U}[\delta]$
subject to the relations $\delta u=\partial^U(u)+(-1)^{|u|}u\delta$ for all $u\in U$.
By definition, this yields formulas like the following:
\begin{gather*}
u(u'\delta)=uu'\delta \\
(u\delta)u'=u\partial^U(u')+(-1)^{|u'|}uu'\delta\\
(u\delta)(u'\delta)=u\partial^U(u')\delta+(-1)^{|u'|}uu'\delta^2=u\partial^U(u')\delta
\end{gather*}
\end{notn}

Next, we summarize some important connections between $U$ and $V$. 
The verifications are straightforward, so we omit them.

\begin{fact}\label{fact161119a}
\begin{enumerate}[(a)]
\item \label{fact161119a1}
$V$ is an associative algebra concentrated in degrees $-1,0,\ldots,q$.
Note that the relation $\delta u=\partial^U(u)+(-1)^{|u|}u\delta$ shows that $V$ is not
in general graded commutative. 
\item \label{fact161119a2}
Given a graded $F$-vector space $X$,
the graded left $V$-modules with underlying $F$-vector space $X$ are in bijection with
the DG $U$-modules with underlying $F$-vector space $X$.
Specifically, given a DG $U$-module $M$, one defines a graded $V$-module structure on it
via the formula $\delta m=\partial^M(m)$.
Conversely, given a graded $V$-module $N$, one defines a DG $U$-module structure on it
via the same formula: $\partial^N(n)=\delta n$.
\end{enumerate}
\end{fact}

\begin{notn}\label{notn161119b}
Given a DG $U$-module $M$, we write $\sv{M}$ for the corresponding graded $V$-module. 
We specify the DG $U$-module structure on $M$ by writing $\su{M}$.
Similarly, given a graded $V$-module $N$, we write $\su{N}$ for the corresponding DG $U$-module,
and we specify the graded $V$-module structure on $N$ by writing $\sv{N}$.
\end{notn}

\begin{para}[Proof of Theorem~\ref{tangent space and Ext}] \label{proof161213a}
For each $i\in\bbz$, we have a natural monomorphism
$$
\grHom[V]{\sv M}{\sv M}_i\into\Hom[U]{\su M}{\su M}_i
$$
the image of which is exactly the set of degree-$i$ cycles in $\Hom[U]{\su M}{\su M}$.
The point is that a graded map $f\colon M\to M$ that is  $V$-linear  is also $U$-linear and, moreover,
it is $V$-linear if and only if it is $U$-linear and it is a cycle in the Hom-complex over $U$.
In particular, the case $i=0$ implies that 
the graded $V$-module homomorphisms $\sv M\to\sv M$ are in bijection with the
DG $U$-module morphisms $\su M\to \su M$.
From this
we obtain an abelian group isomorphism
\begin{equation}\label{eq161213a}
\yext[U]1{\su M}{\su N}\cong\gryexta[V]1{\sv M}{\sv N}.
\end{equation}
Indeed, given an exact sequence
\begin{equation}
0\to N\to X\to M\to 0\label{eq161119a}
\end{equation}
of graded $F$-vector spaces, we have the following:
\begin{enumerate}[(1)]
\item \label{fact161119b2a}
the sequence~\eqref{eq161119a} consists of DG $U$-module morphisms if and only if it consists of
$V$-linear maps;
\item \label{fact161119b2b}
an $F$-linear isomorphism between the extension~\eqref{eq161119a} and the trivial extension
$0\to N\to M\oplus N\to M\to 0$ will consist of DG $U$-module morphisms if and only if it consists of
$V$-linear maps, that is, the  extension~\eqref{eq161119a} will be trivial in 
$\yext[U]1{\su M}{\su N}$ if and only if it is trivial in $\gryexta[V]1{\sv M}{\sv N}$;
and
\item \label{fact161119b2c}
the Baer sum in the DG setting over $U$ corresponds directly to the Baer sum in the graded setting over $V$.
\end{enumerate}
The correspondence in Fact~\ref{fact161119a}\eqref{fact161119a2} between DG $U$-modules and graded $V$-modules 
induces an isomorphism of varieties $\od^U(W)\xra\cong\grod^V(W)$.
It is straightforward to show that this is equivariant with respect to the actions of $\ggl(W)_0$.
In particular, this  induces an isomorphism of orbits
$\ggl(W)_0\cdot\su M\xra\cong\ggl(W)_0\cdot\sv M$.
The first of these isomorphisms 
induces an isomorphism of tangent spaces
$\Tamod_{\su M}\xra\cong\Tagrmod_{\sv M}$
which maps the subspace
$\Taglum_{\su M}\subseteq\Tamod_{\su M}$ isomorphically onto the subspace
$\Taglvm_{\sv M}\subseteq\Tagrmod_{\sv M}$.
Thus, we have 
the first isomorphism in the next sequence
\begin{align*}
\Tamod_{\su M}/\Taglum_{\su M}
&\cong\Tagrmod_{\sv M}/\Taglvm_{\sv M}\\
&\cong\gryexta[V]1{\sv M}{\sv M}\\
&\cong\yext[U]1{\su M}{\su M}.
\end{align*}
The other isomorphisms are from
Theorem~\ref{grVoigt}
and the display~\eqref{eq161213a}.
\hfill$\Box$
\end{para}

\section{Answering Vasconcelos' Question}
\label{sec111227b}

In this section we prove Theorem~\ref{cor111115a} from the introduction; see~\ref{proof111227c}.
We also verify a semi-local version of this result in Theorem~\ref{thm111212a}.
We begin with 
a consequence
of the work from 
Section~\ref{sec161119b},
motivated by
Happel's result~\cite{happel:sm} discussed in the introduction. Recall that $s$ and $W$ are
fixed in Notation~\ref{notn160311a};
other notation comes from Appendix~\ref{sec110211a}.

\begin{prop}\label{prop111115a}
We work with the notation from Section~\ref{sec161119b}.
Let $\s_W(U)$ denote the set of isomorphism classes in $\catd(A)$ of
degree-wise finite semi-free semidualizing DG $U$-modules $C$ such that
$s\geq\sup(C)$, $C_i=0$ for all $i<0$, and $\und{(\tau(C)_{(\leq s)})}\cong W$.
Then $\s_W(U)$ is a finite set.
\end{prop}

\begin{proof}
Let $[C],[C']\in\s_W(U)$, and set $M=\tau(C)_{(\leq s)}$ and $M'=\tau(C')_{(\leq s)}$.
Since we have $s\geq\sup(C)$, we conclude that $M\simeq C$ in $\catd(U)$,
and similarly for $M'$.

We  observe that the semidualizing assumption, in particular the condition 
$$\Ext[U]1{M}{M}\cong\Ext[U]1CC=0$$
implies that the orbit $\ggl(W)_0\cdot M$ is open in $\od^U(W)$.
Indeed, from~\cite[Proposition~4.4]{nasseh:egdgm} we have
$\yext[U]1MM=0$, so 
Theorem~\ref{tangent space and Ext} shows that
$\Tamod_M=\Taglm_M$. 
Lemma~\ref{prop160312a}\eqref{prop160312a2} implies that 
$\ggl(W)_0\cdot M$ is locally closed and non-singular.
So, the orbit $\ggl(W)_0\cdot M$ is open by~\cite[Lemma~2.14]{brion:rq}.

Next, note that if $\ggl(W)_0\cdot M=\ggl(W)_0\cdot M'$, then $[C]=[C']$:
indeed, by Proposition~\ref{prop160311a}, we have
$C\simeq M\cong M'\simeq C'$.

As $\od^U(W)$ is quasi-compact, it can only have finitely many open orbits.
By what we have already shown, the map $[C]\mapsto \ggl(W)\cdot M$ is a well-defined
injective map from $\s_W(U)$ into a finite set of open orbits, so $\s_W(U)$ is finite, as desired.
\end{proof}

\begin{para}[Proof of Theorem~\ref{cor111115a}] \label{proof111227c}
A result of Grothendieck~\cite[Proposition (0.10.3.1)]{grothendieck:ega3-1} provides
a flat local ring homomorphism $(R,\m,k)\to(R',\m',k')$ such that
$k'$ is algebraically closed. Composing with the natural map from $R'$ to its
$\m'$-adic completion, we assume that $R'$ is complete.
By~\cite[Theorem II(c)]{frankild:rrhffd},  the induced map
$\s(R)\to\s(R')$ is injective.
Thus it suffices to prove the result for $R'$, so   assume
that $R$ is complete with algebraically closed residue field.

Let $\underline{t}=t_1,\cdots,t_n$ be a minimal generating sequence
for $\mathfrak m$,
and set $K=K^R(\underline t)$, the Koszul complex.
The  map
$\s(R)\to\s(K)$ induced by $C\mapsto\Otimes KC$ is  bijective 
by~\cite[Corollary 3.10]{nasseh:ldgm}.
Thus, it suffices to show that $\s(K)$ is finite.
Note that for each semidualizing $R$-complex $C$, we have
$\amp(C)\leq\dim(R)-\depth(R)$ by~\cite[(3.4) Corollary]{christensen:scatac}.
A standard result about $K$ (see, e.g., \cite[1.3]{foxby:daafuc}) implies that
\begin{equation}\label{eq111201a}
\amp(\Otimes KC)\leq\amp(C)+n\leq\dim(R)-\depth(R)+n.
\end{equation}
Set $s=\dim(R)-\depth(R)+n$.

From~\cite[(2.8)]{avramov:lgh} there is a  finite dimensional DG $k$-algebra $U$, as in Notation~\ref{notn160311a} with $F=k$, that is
linked to $K$ by a sequence of (quasi)isomorphisms of 
local
DG algebras. 
By Lemma~\ref{lem111117a}\eqref{lem111117a4} there is a bijection between $\s(K)$ and $\s(U)$.
Thus, it suffices to show that $\s(U)$ is finite.

Let $C'$ be a semidualizing DG $U$-module,
and let $C$ be a semidualizing $R$-complex
corresponding to $C'$ under the bijection
given above. 
From Lemma~\ref{lem111117a}\eqref{lem111117a1}
and the display~\eqref{eq111201a}, we have 
$\amp(C')=\amp(\Otimes KC)\leq s$.
By applying an appropriate shift
we assume without loss of generality that $\inf(C)=0=\inf(C')$,
so we have $\sup(C')\leq s$. 
Let $L\xra\simeq C'$ be a minimal semi-free  resolution of $C'$ over $U$;
see Fact~\ref{fact110218d}.
The conditions $\sup(L)=\sup(C')\leq s$ imply that $L$ (and hence $C'$)
is quasiisomorphic to the truncation $\wti L:=\tau(L)_{(\leq s)}$.
We set $W:=\und{\wti{L}}$
and work with the notation set in Section~\ref{sec161119b}.

We claim that there is an integer $\lambda\geq 0$, depending only on $R$
and $U$, such that 
$\sum_{i=0}^sr_i\leq \lambda$.
(Recall that $r_i$ and other quantities are fixed in Notation~\ref{notn160311a}.)
To see this, first note that for $i=0,
\ldots,s$
we have $L_i=\bigoplus_{j=0}^i U^{\beta_{i-j}^U(C')}_{j}$;
see Fact~\ref{fact110218d}.
Using minimal semi-free resolutions and Lemma~\ref{lem131221a}, we conclude that
$$\beta_{j}^U(C')=\beta_{j}^R(C)\leq\mu^{j+\depth(R)}_R(R)$$
for all $j$.
With our notation $n_i=\dim_F(U_i)$ from~\ref{notn160311a},
it follows that 
\begin{equation*}
r_i
=\rank_F(\widetilde{L}_i)
\leq\rank_F(L_i)
=\sum_{j=0}^i n_{i-j}\beta_{j}^U(C')
\leq\sum_{j=0}^i n_{i-j}\mu^{j+\depth(R)}_R(R).
\end{equation*}
And we conclude that
\begin{equation*}
\sum_{i=0}^sr_i
\leq \sum_{i=0}^s\sum_{j=0}^i n_{i-j}\mu^{j+\depth(R)}_R(R).
\end{equation*}
Since the numbers in the right hand side of this inequality only depend on $R$ and $U$,
we have found the desired value for $\lambda$.

Because there are only finitely many $(r_0,\ldots,r_s)\in \mathbb{N}^{s+1}$ with 
$\sum_{i=0}^sr_i\leq \lambda$,
there are only finitely many $W$ that occur from this construction,
say $W^{(1)},\ldots,W^{(b)}$.
Proposition~\ref{prop111115a}
implies that $\s(U)=\s_{W^{(1)}}(U)\cup\cdots\cup \s_{W^{(b)}}(U)$
is finite.
\hfill$\Box$
\end{para}

\begin{disc}\label{disc131220a}
Note that Theorem~\ref{cor111115a} provides a positive answer to Question~\ref{q111227a}
since one has $\s_0(R)\subseteq\s(R)$; 
see Definition~\ref{dfn:DGsdm}.
\end{disc}

Now, we move toward our semi-local version of Theorem~\ref{cor111115a}.

\begin{defn}
The \emph{non-Gorenstein locus} of $R$ is
$$\nGor(R):=\{\text{maximal ideals $\m\subset R$}\mid\text{$R_{\m}$ is not Gorenstein}\}\subseteq\mspec(R)$$
where $\mspec(R)$ is the set of maximal ideals of $R$.
\end{defn}

\begin{disc}\label{disc111226a}
For ``nice'' rings, e.g. rings with a dualizing complex~\cite{sharp:d}, 
the set $\nGor(R)$ is closed
in $\mspec(R)$.
\end{disc}

\begin{thm}\label{thm111212a}
Assume that $R$ satisfies one of the following conditions:
\begin{enumerate}[\rm(1)]
\item \label{thm111212a1} $R$ is semi-local, or
\item \label{thm111212a2} $R$ is Cohen-Macaulay and
$\nGor(R)$ is finite.
\end{enumerate}
Then the sets $\ol{\s_0}(R)$ and $\ol{\s}(R)$ are finite; 
see Notation~\ref{notn111226a}.
\end{thm}

\begin{proof}
Because of the containment $\ol{\s_0}(R)\subseteq\ol{\s}(R)$, it suffices to show that
$\ol{\s}(R)$ is finite.
Given
a finite subset
$X=\{\m_1,\ldots,\m_n\}\subseteq\mspec(R)$, let
$f_X
\colon \ol{\s}(R)\to\prod_{i=1}^n\s(R_{\m_i})$
be given by the formula $f_X
(\orbit{C}):=([C_{\m_1}],\ldots,[C_{\m_n}])$.
This is well-defined 
by Fact~\ref{fact111227a}.

In each case~\eqref{thm111212a1}--\eqref{thm111212a2} 
we  show that there is a finite set $X$
such that $f_X$ is injective.
Then Theorem~\ref{cor111115a} implies that 
the set $\prod_{i=1}^n\s(R_{\m_i})$ is finite, so
$\ol{\s}(R)$ is also finite.

\eqref{thm111212a1}
Assume that $R$ is semi-local, and set $X:=\mspec(R)$.
To show that 
$f_X$ is injective,
let $[B],[C]\in\s(R)$ be
such  that  $B_{\m_i}\sim C_{\m_i}$ for $i=1,\ldots,n$.
According to Fact~\ref{fact111227a}, to show that $\orbit C=\orbit B$,
it suffices to show that $\Ext jBC=0$ for $j\gg 0$. 
Since we know that $\Ext[R_{\m_i}]{j}{B_{\m_i}}{B_{\m_i}}=0$ for all $j\geq 1$,
we conclude that there are integers $j_1,\ldots,j_n$ such that
for $i=1,\ldots,n$ we have
$\Ext[R_{\m_i}]{j}{B_{\m_i}}{C_{\m_i}}=0$ for all $j\geq j_i$.
Since $B$ is homologically finite, we have
$0=\Ext[R_{\m_i}]{j}{B_{\m_i}}{C_{\m_i}}\cong\Ext jBC_{\m_i}$
for all $j\geq \max_i{j_i}$. Since vanishing is a local property, it follows that
$\Ext jBC=0$ for $j\gg 0$, as desired.

\eqref{thm111212a2}
Now, assume that $R$ is Cohen-Macaulay and
$\nGor(R)$ is finite, and set
$X:=\nGor(R)$.
To show that $f_X$ is injective,
let $[B],[C]\in\s(R)$ be
such  that $B_{\m_i}\sim C_{\m_i}$ for $i=1,\ldots,n$.
Lemma~\ref{lem111227a} provides
semidualizing modules $B'$ and $C'$
such that  $\orbit{B}=\orbit{B'}$
and $\orbit{C}=\orbit{C'}$,
so we assume without loss of generality
that $B$ and $C$ are modules.

We claim that 
$B_{\m}\sim C_{\m}$ for all maximal ideals $\m\subset R$ and $\Ext iBC=0$ for all $i\geq 1$.
(Then the desired conclusion then
follows 
from Fact~\ref{fact111227a}.)
As
$B$ is a finitely generated $R$-module, it suffices to show that
$B_{\m}\sim C_{\m}$  and 
$\Ext[R_{\m}] i{B_{\m}}{C_{\m}}=0$ for all $i\geq 1$ and for all maximal ideals $\m\subset R$.

Case 1:  $\m\in\nGor(R)$. In this case, we have $B_{\m}\sim C_{\m}$,
by assumption. Since $B$ and $C$ are both modules, this implies that
$B_{\m}\cong C_{\m}$, so the fact that $B_{\m}$ is semidualizing
over $R_{\m}$ implies that
$\Ext[R_{\m}] i{B_{\m}}{C_{\m}}\cong \Ext[R_{\m}] i{B_{\m}}{B_{\m}}=0$.

Case 2: $\m\notin\nGor(R)$. In this case, the ring $R_{\m}$ is 
Gorenstein, so we have $B_{\m}\cong R_{\m}\cong C_{\m}$ 
by~\cite[(8.6) Corollary]{christensen:scatac},
and the desired vanishing follows.
\end{proof}

\appendix
\section{DG Modules}
\label{sec110211a}

We assume that the reader is familiar with the category of $R$-complexes.
For clarity, though, we include some  notation.

\begin{defn}
  \label{cx}
We index chain complexes of $R$-modules (``$R$-complexes'' for short)  homologically:
$M = \cdots \xra{\partial_{n+2}^{M}} M_{n+1} \xra{\partial_{n+1}^{M}} M_n
\xra{\partial_{n}^{M}} M_{n-1} \xra{\partial_{n-1}^{M}} \cdots$.
The degree of an element $m\in M$ is denoted $|m|$.
The \emph{tensor product} of two $R$-complexes $M,N$ is denoted $\Otimes MN$,
and the \emph{Hom complex} is denoted $\Hom MN$.
A \emph{chain map} $M\to N$ is a cycle in $\Hom MN_0$.
\end{defn}

Next we discuss  DG algebras, which are treated  in, e.g., 
\cite{avramov:ifr,avramov:dgha,avramov:bsolrhoffd,beck:sgidgca,felix:rht,keller:ddgc}.
We follow the notation and terminology from~\cite{avramov:dgha,beck:sgidgca}; given the slight differences in the literature,
though, we include  a summary next.

\begin{defn}
  \label{DGK}
A \emph{positively graded commutative differential graded $R$-algebra} (\emph{DG $R$-algebra} for short)
is an $R$-complex $A$ equipped with a
chain map
$\mu^A\colon A\otimes_RA\to A$ 
such that the product
$ab:=\mu^A(a\otimes b)$
is associative, unital, and graded commutative, 
and
such that $A_i=0$ for $i<0$.
The map $\mu^A$ is the \emph{product} on $A$.
Given a DG $R$-algebra $A$, the \emph{underlying algebra} is the
graded commutative  $R$-algebra
$\und{A}=\oplus_{i=0}^\infty A_i$.
When $R$ is a field and $\rank_R(\und A)<\infty$, then $A$ is \emph{finite dimensional} over $R$.
We say that $A$ is \emph{homologically degree-wise noetherian}
if $\HH_0(A)$ is noetherian and the $\HH_0(A)$-module $\HH_i(A)$ is  finitely generated for all $i\geq 0$.

A \emph{morphism} of DG $R$-algebras is a chain map
$f\colon A\to B$ between DG $R$-algebras  respecting products and multiplicative identities.
A \emph{quasiisomorphism} of DG $R$-algebras is a morphism that is  a quasiisomorphism.

Assume that $(R,\m)$ is local.
We say that $A$ is  \emph{local} if it is homologically degree-wise noetherian and
the ring $\HH_0(A)$ is a local $R$-algebra, so,
$\HH_0(A)$ is a local  ring 
with
maximal ideal $\m_{\HH_0(A)}$ 
containing
$\m\HH_0(A)$.
In this case, the
composition $A\to \HH_0(A)\to \HH_0(A)/\m_{\HH_0(A)}$
is a surjective morphism of DG $R$-algebras with kernel of the form
$\m_A=\cdots \xra{\partial_2^A} A_1 \xra{\partial_1^A} \mathfrak m_0\to 0$
where $\m_0$ is the preimage of $\m_{\HH_0(A)}$ in $A_0$. 
The
ideal $\m_0$ is maximal and we have $\m A_0\subseteq\m_0$.
The DG ideal $\m_A$  is the \emph{augmentation ideal} of $A$.
\end{defn}

For this paper, an important example is the next one.

\begin{ex}
Given a sequence $\ba=a_1,\cdots,a_n\in R$,
the Koszul complex $K=K^R(\ba)$ is a DG $R$-algebra with product given by the exterior (i.e., wedge) product.
If $(R,\m)$ is local
and $\ba\in\m$, then $K$ is a local
DG $R$-algebra
with augmentation ideal
$\m_K=(0\to R\to\cdots\to R^n\to\m\to 0)$.
\end{ex}

In the passage to DG algebras, we  
focus
on DG modules,
described next.

\begin{defn}
Let $A$ be a DG $R$-algebra. 
A 
\emph{DG $A$-module} is
an $R$-complex $M$ with a
chain map
$\mu^M\colon \Otimes AM\to M$ such that the rule $am:=\mu^M(a\otimes m)$
is associative and unital.
The map $\mu^M$ is the \emph{scalar multiplication} on $M$.
The \emph{underlying $\und{A}$-module} associated to $M$ is the
$\und{A}$-module
$\und{M}=\oplus_{i\in\bbz} M_i$.

Given  DG $A$-modules $M$, $N$,
the DG $A$-modules $\Otimes[A]MN$ and $\Hom[A]MN$ are defined  in~\cite[Section~1]{avramov:bsolrhoffd}.
A \emph{morphism} of DG $A$-modules is a cycle in $\Hom[A]MN_0$.
Isomorphisms in the category of DG $A$-modules are identified by the
symbol $\cong$. 
Quasiisomorphisms of DG $A$-modules are identified by the
symbol $\simeq$.
\end{defn}

\begin{ex}\label{ex110218a}
Consider the ring $R$ as a DG $R$-algebra.
A DG $R$-module is just an $R$-complex, and a morphism of DG $R$-modules is simply a chain map.
\end{ex}

\begin{defn}
  \label{DGK2}
Let $A$ be a DG $R$-algebra,  let $i$ be an integer, and let 
$M$ be a DG $A$-module.
The $i$th \emph{suspension} of  $M$
is the DG $A$-module $\shift^iM$ defined by $(\shift^iM)_n :=
M_{n-i}$ and $\partial^{\shift^iM}_n := (-1)^i\partial^M_{n-i}$. 
The scalar multiplication
on $\shift^iM$ is defined as
$\mu^{\shift^iM}(a\otimes m):=(-1)^{i|a|}\mu^M(a\otimes m)$.
\end{defn}

\begin{defn}\label{disc110616a}
Let $A$ be a    DG $R$-algebra.
The category of DG $A$-modules  described above
is an  abelian category; see, e.g., 
\cite[Introduction]{keller:ddgc}.
So, given DG $A$-modules $M$ and $N$, the \emph{Yoneda Ext group} $\yext[A] 1MN$,
defined as the set of equivalence classes of 
exact sequences 
$0\to N\to X\to M\to 0$
of DG $A$-modules,
is a well-defined abelian group under the Baer sum; see, e.g., \cite[(3.4.6)]{weibel:iha}.
An explicit description of the Baer sum in this setting is given in the proof of~\cite[Theorem~3.5]{nasseh:egdgm}.
\end{defn}

In the remainder of this appendix, we provide some background on semidualizing DG-modules, e.g., semidualizing complexes.
We assume that the reader is familiar with 
the derived category $\catd(R)$.
References
for this
include~\cite{christensen:gd,christensen:dcmca,gelfand:moha,hartshorne:rad,verdier:cd,verdier:1}.
For clarity, we include some definitions and notation.

\begin{defn}
  \label{cx'}
Let $A$ be a DG $R$-algebra.
The \emph{infimum}, \emph{supremum}, and \emph{amplitude} of a DG $A$-module
$M$ are
\begin{align*}
\inf(M)
&:=\inf\{n\in\bbz\mid\HH_n(M)\neq 0\}
\\
\sup(M)
&:=\sup\{n\in\bbz\mid\HH_n(M)\neq 0\}\\
\amp(M)
&:=\sup(M)-\inf(M).
\end{align*}
The
DG $A$-module
$M$ is \emph{bounded below} if $M_n=0$ for all $n\ll 0$;
it is \emph{degree-wise finite} if $M_i$ is finitely generated over $A_0$ for each $i$;
it  is \emph{homologically bounded} if $\HH_i(M)=0$ for $|i|\gg 0$;
it is \emph{homologically degree-wise finite} if each $\HH_0(A)$-module
$\HH_n(M)$ is finitely generated; and it is \emph{homologically
finite} if it is homologically both bounded and degree-wise finite.

A DG $A$-module $Q$ is 
\emph{semi-projective} if $\Hom[A]Q-$ respects surjective quasiisomorphisms,
that is, if $\und Q$ is a projective graded $\und R$-module and $\Hom[A]Q-$ respects quasiisomorphisms.
A DG $A$-module $L$ is \emph{semi-free} if the graded $A^\natural$-module $L^ \natural$ 
has a graded basis $X$ which decomposes as a disjoint union $X=\sqcup_{i=0}^\infty X^i$ such that
$\partial^L(X^i)$ is contained in the DG submodule $RX^{i-1}$ for each $i\geq 0$, where $X^{-1}:=\emptyset$.
A \emph{semi-free (resp., semi-projective) resolution} of  $M$ is a quasiisomorphism
$L\xra{\simeq}M$ of DG $A$-modules such that $L$ is semi-free (resp., semi-projective).
If $R$ and $A$ are local, then
a \emph{minimal semi-free resolution} of $M$ is
a semi-free resolution $L\xra\simeq M$ such that 
each  graded basis of $L^\natural$ is finite in each degree and
the differential on
$\Otimes[A]{(A/\m_A)}{L}$ is 0.
\end{defn}

\begin{defn}
  \label{DGK2'}
Let $A$ be a DG $R$-algebra.
The derived category $\catd(A)$ is formed from the category of DG $A$-modules
by formally inverting the quasiisomorphisms; see~\cite{keller:ddgc}.
Isomorphisms in $\catd(A)$  are identified by the symbol $\simeq$,
and isomorphisms up to shift in $\catd(A)$  are identified by  $\sim$.

The  derived functors $\Lotimes[A]MN$ and $\Rhom[A]MN$ are given, e.g., via a semi-projective
resolution $P\xra\simeq M$,
as $\Lotimes[A]MN\simeq\Otimes[A]PN$ and $\Rhom[A]MN\simeq\Hom[A]PN$.
For each $i\in\bbz$,
set $\Tor[A]iMN:=\HH_i(\Lotimes[A]MN)$ and $\Ext[A]iMN:=\HH_{-i}(\Rhom[A]MN)$.\footnote{See~\cite{nasseh:egdgm} for 
\label{footnote161213a}
differences and similarities between,
e.g.,
$\Ext[A]1MN$ and 
$\yext[A] 1MN$.}
\end{defn}

\begin{fact}\label{fact110218d}
Let $A$ be a DG $R$-algebra.
If $L$ is a bounded below DG $A$-module such that $L^\natural$ is a graded free $A^\natural$-module,
then $L$ is semi-free.
If $M$ is a DG $A$-module such that $j=\inf(M)>-\infty$, then $M$ has  a semi-free (hence, semi-projective) resolution
$L\xra\simeq M$ such that $L_i=0$ for all $i<j$.
If $A$ is homologically degree-wise noetherian and $M$ is homologically degree-wise finite, then $L$ can be chosen so that
$\und{L}\cong\oplus_{i=j}^\infty \shift^i(\und{A})^{\beta_i}$
for some integers $\beta_i$.
If, in addition,
$R$ 
and $A$ are local, then $L$ can be chosen to be  minimal by~\cite[Lemma A.3.iii]{felix:gs}.
\end{fact}

In the passage from $R$ to $U$ in our proof of Theorem~\ref{cor111115a},
we use Christensen and Sather-Wagstaff's notion of semidualizing DG $U$-modules
from~\cite{christensen:dvke},
defined next.

\begin{defn}
\label{dfn:DGsdm}
Let $A$ be a DG $R$-algebra, and let $M$ be a DG $A$-module.
The \emph{homothety morphism} $X^A_M\colon A\to \HomA MM$ is given by
$X^A_M(a)(m)=am$.
This induces a \emph{homothety morphism} $\chi^A_M\colon A\to \RhomA MM$ in $\catd(A)$.

Assume that $A$ is homologically degree-wise noetherian. Then  $M$ is a \emph{semidualizing}
DG $A$-module if $M$ is homologically finite
and the homothety morphism
$\chi^A_M\colon A\to\Rhom[A]{M}{M}$ is an isomorphism in $\catd(A)$.\footnote{Note here that the existence of a
semidualizing DG $A$-module forces $A$ to be homologically bounded. Conversely, if $A$ is homologically bounded,
then it has a semidualizing DG $A$-module, namely, $A$ itself. See, e.g., \cite{altmann:sdmtp}.}
Let $\s(A)$ denote the set of shift-isomorphism classes 
in $\catd(A)$ of semidualizing DG $A$-modules,
that is, the set of equivalence classes of semidualizing DG $A$-modules 
under the relation $\sim$ from Definition~\ref{DGK2'}.
Over the ring $R$, semidualizing DG $R$-modules are called \emph{semidualizing complexes};
in this setting, we let $\s_0(R)$ denote the
set of  isomorphism classes of semidualizing $R$-modules,
which is naturally a subset of $\s(R)$.
\end{defn}

The following base-change 
construction is
used 
in the passage from $R$ to $U$ in our proof of Theorem~\ref{cor111115a}.

\begin{disc}\label{disc111223a}
Let $A\to B$ be a morphism of DG $R$-algebras, and let $M$ and $N$ be 
DG $A$-modules.
The
``base changed'' complex $\Otimes[A]BM$ has the structure of a DG $B$-module
by the action $b(\Otimes[] {b'}m):=\Otimes[]{(bb')}{m}$.
This structure is compatible with the DG $A$-module structure
on $\Otimes[A]BM$ via restriction of scalars.
Furthermore, this induces a well-defined operation
$\catd(A)\to\catd(B)$ given by
$M\mapsto\Lotimes[A]BM$.

Given an element $f\in\Hom[A]MN_i$, define
$\Otimes[A]Bf\in \Hom[B]{\Otimes[A]BM}{\Otimes[A]BN}_i$ by the formula
$(\Otimes[A]Bf)(\Otimes[]bm):=(-1)^{i|b|}\Otimes[]b{f(m)}$.
This yields a morphism of DG $A$-modules 
$\Hom[A]MN\to \Hom[B]{\Otimes[A]BM}{\Otimes[A]BN}$
given by $f\mapsto \Otimes[A]Bf$ which, in turn, provides a  morphism 
$\Rhom[A]MN\to  \Rhom[B]{\Lotimes[A]BM}{\Lotimes[A]BN}$
in $\catd(A)$.
\end{disc}

The next lemma is essentially from~\cite{keller:ilchdga} and~\cite{pauksztello:hedga}.

\begin{lem}\label{lem111117a}
Let $\vf\colon A\to B$ be a quasiisomorphism of DG $R$-algebras.
\begin{enumerate}[\rm(a)]
\item \label{lem111117a3}
The base change functor $\Lotimes[A]B-$ induces an equivalence of derived
categories $\catd(A)\to\catd(B)$ whose quasi-inverse is given by restriction of scalars.
\item \label{lem111117a1}
For each DG $A$-module $X\in\catd(A)$, one has
$X\simeq \Lotimes[A]BX$ in $\catd(A)$.
\item \label{lem111117a4}
The equivalence from part~\eqref{lem111117a3} induces a bijection from $\s(A)$ to $\s(B)$.
\end{enumerate}
\end{lem}

\begin{proof}
\eqref{lem111117a3}
See, e.g., \cite[7.6 Example]{keller:ilchdga}.

\eqref{lem111117a1}
The equivalence from part~\eqref{lem111117a3}
implies that the natural morphism $X\to\Lotimes[A]BX$ is an isomorphism in $\catd(A)$.

\eqref{lem111117a4}
Let $X$ be a  DG $A$-module. 
We show that $X$ is a semidualizing DG $A$-module if and only if
$\Lotimes[A]BX$ is a semidualizing DG $B$-module.
As the maps $X\to\Lotimes[A]BX$ and $A\to B$ are (quasi)isomorphisms,
it follows that $X$ is homologially finite over $A$ if and only if 
$\Lotimes[A]BX$ is homologially finite over $B$.
It remains to show that the  morphism
$\chi^A_X\colon A\to\Rhom[A]{X}{X}$ is an isomorphism in $\catd(A)$
if and only if 
$\chi^B_{\Lotimes[A]BX}\colon B\to\Rhom[B]{\Lotimes[A]BX}{\Lotimes[A]BX}$ is an isomorphism
in $\catd(B)$.
It is routine to show that 
the following diagram commutes
$$\xymatrix{
A\ar[r]^-{\chi^A_X}\ar[d]_{\vf}^\simeq
&\Rhom[A]{X}{X}\ar[d]^{\omega}_\simeq \\
B\ar[r]_-{\chi^B_{\Lotimes[A]BX}}
&\Rhom[B]{\Lotimes[A]BX}{\Lotimes[A]BX}}$$
where  $\omega$ is the morphism
from Remark~\ref{disc111223a}.
As $\omega$ is an isomorphism by~\cite[Proposition 2.1]{pauksztello:hedga},
the desired equivalence follows.
\end{proof}

\begin{defn}\label{truncations}
Let $A$ be a  DG $R$-algebra, and let $M$ be a DG $A$-module.
Given an integer $n$, 
the \emph{$n$th soft left truncation of $M$} is the complex
$$
\tau(M)_{(\leq n)}:=\cdots\to 0\to M_n/\im(\partial^M_{n+1})\to M_{n-1}\to M_{n-2}\to \cdots
$$
with differential induced by $\partial^M$.
In other words, $\tau(M)_{(\leq n)}$ is the quotient DG $A$-module
$M/M'$ where $M'$ is the following DG submodule of $M$:
$$M'=\cdots\to M_{n+2}\to M_{n+1}\to\im(\partial^M_{n+1})\to 0.$$
Note that $M'\simeq 0$ if and only if $n\geq \sup(M)$,
so the natural morphism $M\to \tau(M)_{(\leq n)}$ of DG $A$-modules
is a quasiisomorphism if and only if $n\geq \sup(M)$.
\end{defn}

\begin{defn}\label{defn111226a}
Let $A$ be a local DG $R$-algebra with $k=A/\m_A$, and let $M$ be a homologically finite DG $A$-module.
For each integer $i$, the $i$th 
\emph{Betti number} and the $i$th \emph{Bass number} are respectively
\begin{align*}
\beta_{i}^A(M)&:=\rank_k(\Tor[A]{i}{k}{M})&
\mu^{i}_A(M)&:=\rank_k(\Ext[A]{i}{k}{M}).
\intertext{The \emph{Poincar\'e series} and the \emph{Bass series} of $M$ are the formal Laurent series}
P^M_A(t)&:=\sum_{i\in\bbz}\beta_{i}^A(M)t^i&
I_M^A(t)&:=\sum_{i\in\bbz}\mu^{i}_A(M)t^i.
\end{align*}
\end{defn}

The inequality in the next lemma may be of independent interest.

\begin{lem}\label{lem131221a}
Assume that $R$ is local, and let $C$ be a semidualizing $R$-complex such that $\inf(C)=0$.
Then one has $\beta^R_p(C)\leq \mu^{p+\depth R}_R(R)$ for all $p\geq 0$.
\end{lem}

\begin{proof}
The isomorphism $\Rhom[R]{C}{C}\simeq R$ implies the following
equality of power series
$I_R^R(t)=P_R^C(t)I_C^R(t)$.
See~\cite[(1.5.3)]{avramov:rhafgd}.
We conclude that for each $m$ we have 
$$
\mu^m_R(R)=\sum_{i=0}^m \beta^R_i(C)\mu_R^{m-i}(C).
$$
In particular, for $m<\depth(R)$, we have
$$
0=\mu^{m}_R(R)
\geq \beta^R_0(C)\mu_R^{m}(C).
$$
The equality $\inf(C)=0$ implies that $\beta_0^R(C)\neq 0$
by~\cite[(1.7.1)]{christensen:scatac}, 
so it follows that $\mu_R^{m}(C)=0$.
For $m=\depth(R)$, we conclude from this that
$$
0\neq\mu^{\depth(R)}_R(R)
= \beta^R_0(C)\mu_R^{\depth(R)}(C)
$$
and hence $\mu_R^{\depth(R)}(C)\neq 0$.
Similarly, for $m=p+\depth(R)$, we have
\begin{align*}
\mu^{p+\depth(R)}_R(R)
&\geq \beta^R_p(C)\mu_R^{\depth(R)}(C)
\geq \beta^R_p(C)
\end{align*}
as desired.
\end{proof}

We note that, over a non-local ring, the set $\s(R)$ may not be finite.
For instance, the \emph{Picard group} $\Pic(R)$, consisting of finitely generated
rank-1 projective $R$-modules, is contained in
$\s_0(R)\subseteq\s(R)$, so $\s(R)$ can  be infinite even
when $R$
is a  Dedekind domain. We use some  notions from~\cite{frankild:rbsc}
to deal with this.

\begin{defn}\label{defn}
A \emph{tilting $R$-complex} is a semidualizing $R$-complex of 
finite projective dimension.
The \emph{derived Picard group of $R$} is the set $\dpic(R)$ of isomoprhism
classes in $\catd(R)$ of tilting $R$-complexes.
The isomorphism class of a tilting $R$-complex $L$ is denoted $[L]\in\dpic(R)$.
\end{defn}

\begin{disc}\label{disc111212a}
A homologically finite $R$-complex $L$ is tilting if and only if $L_{\m}\sim R_{\m}$ for all
maximal (equivalently, for all prime) ideals $\m\subset R$, 
by~\cite[Proposition 4.4 and Remark 4.7]{frankild:rbsc}.
In~\cite{avramov:rrc1} tilting complexes are called ``invertible''.
The derived Picard group $\dpic(R)$ is an abelian group under the operation
$[L][L']:=[\Lotimes L{L'}]$. The identity in $\dpic(R)$ is $[R]$, and $[L]^{-1}=[\Rhom LR]$.
The classical Picard group $\Pic(R)$ is naturally a subgroup of $\dpic(R)$.
The group $\dpic(R)$ acts on $\s(R)$ in a natural way:
$[L][C]:=[\Lotimes L{C}]$.
See~\cite[Properties~4.3 and Remark 4.9]{frankild:rbsc}.
This action restricts to an action of $\Pic(R)$ on $\s_0(R)$ given by
$[L][C]:=[\Otimes L{C}]$.
\end{disc}

\begin{notation}\label{notn111226a}
The set of orbits in $\s(R)$ under the action of $\dpic(R)$ is denoted
$\ol{\s}(R)$,\footnote{Observe that the notations $\s(R)$ and $\ol{\s}(R)$ 
represent different sets in~\cite{frankild:rbsc}.} 
and the set of orbits in $\s_0(R)$ under the action of $\Pic(R)$ is denoted
$\ol{\s_0}(R)$.
\end{notation}

\begin{fact}\label{fact111227a}
Given semidualizing $R$-complexes $A$ and $B$,
the following conditions are equivalent by~\cite[Proposition 5.1]{frankild:rbsc}:
\begin{enumerate}[(i)]
\item the orbits $\orbit A$ and $\orbit B$ are equal, i.e.,
there is an element $[P]\in\dpic(R)$ such that $B\simeq \Lotimes PA$; and
\item $A_{\m}\sim B_{\m}$ for all maximal ideals $\m\subset R$ and 
$\Ext iAB=0$ for $i\gg 0$.
\end{enumerate}
It is straightforward to show that the natural inclusion 
$\s_0(R)\subseteq\s(R)$ gives an inclusion $\ol{\s_0}(R)\subseteq\ol{\s}(R)$.
\end{fact}

\begin{lem}\label{lem111227a}
Assume that $R$ is Cohen-Macaulay (not necessarily local), 
and let $C$ be a semidualizing $R$-complex.
There is an element $[L]\in\dpic(R)$ such that 
$\Lotimes LC$ is isomorphic in $\catd(R)$ to a module.
In other words, the orbit $\orbit C$ contains a module.
\end{lem}

\begin{proof}
For each $\p\in\spec(R)$, the fact that $R$ is Cohen-Macaulay implies
that $\amp(C_{\p})=0$ by~\cite[(3.4) Corollary]{christensen:scatac}, 
that is, $C_{\p}\sim\HH_i(C)_{\p}\neq 0$ for some $i$.
As $\amp(C)<\infty$, this implies that
$\spec(R)$ is the disjoint union 
$$\spec(R)=\bigcup_{i=\inf(C)}^{\sup(C)}\supp_R(\HH_{i}(C)).$$
It follows that each set $\supp_R(\HH_{i}(C))$ is both open and closed.
So, if $\supp_R(\HH_{i}(C))$ is non-empty, then it is a union of connected
components of $\spec(R)$.

Let $e_1,\ldots,e_p$ be a ``complete set of orthogonal primitive idempotents of
$R$'' as in~\cite[4.8]{avramov:rrc1}. Then $R\cong R_{e_1}\times\cdots\times R_{e_p}$
and each $\spec(R_{e_i})$ is naturally homeomorphic to a connected component
of $\spec(R)$. From the previous paragraph, for $i=1,\ldots,p$ we have
$C_{e_i}\simeq\shift^{u_i}\HH_{u_i}(C_{e_i})$, and $\HH_{u_i}(C_{e_i})$
is a semidualizing $R_{e_i}$-module. Each $R$-module $M$
has a natural decomposition $M\cong\oplus_{i=1}^pM_{e_i}$
that is compatible with the product decomposition of $R$,
and it follows that $C\simeq \oplus_{i=1}^p\shift^{u_i}\HH_{u_i}(C_{e_i})$.

Let $L=\oplus_{i=1}^p\shift^{-u_i}R_{e_i}$. 
Then $L$ is a tilting $R$-complex by Remark~\ref{disc111212a}, and
\begin{align*}
\Lotimes LC
&\simeq\Lotimes{(\oplus_{i=1}^p\shift^{-u_i}R_{e_i})}{(\oplus_{i=1}^p\shift^{u_i}\HH_{u_i}(C_{e_i}))}\\
&\simeq\oplus_{i=1}^p
\Lotimes[R_{e_i}]{(\shift^{-u_i}R_{e_i})}{(\shift^{u_i}\HH_{u_i}(C_{e_i}))}\\
&\simeq\oplus_{i=1}^p
\HH_{u_i}(C_{e_i}).
\end{align*}
Since 
$\oplus_{i=1}^p
\HH_{u_i}(C_{e_i})$ is an $R$-module, this establishes the lemma.
\end{proof}

\providecommand{\bysame}{\leavevmode\hbox to3em{\hrulefill}\thinspace}
\providecommand{\MR}{\relax\ifhmode\unskip\space\fi MR }
\providecommand{\MRhref}[2]{%
  \href{http://www.ams.org/mathscinet-getitem?mr=#1}{#2}
}
\providecommand{\href}[2]{#2}

\affiliationone{Saeed Nasseh, 
Department of Mathematical Sciences,
Georgia Southern University,
Statesboro, Georgia 30460, USA}
\email{snasseh@georgiasouthern.edu}

\affiliationtwo{Sean Sather-Wagstaff, 
Department of Mathematical Sciences, 
Clemson University, 
O-110 Martin Hall, Box 340975, 
Clemson, S.C. 29634, 
USA}
\email{ssather@clemson.edu}

\end{document}